\documentclass{article}
\usepackage{latexsym}
\usepackage{amsmath}
\usepackage{amssymb}
\usepackage[dvips]{graphicx}

\newcommand{\rtre}{{\bf R}^3}
\newcommand{\rdue}{{\bf R}^2}
\newcommand{\erre}{{\bf R}}

\newcommand{\graph}{{\rm graph}}
\newcommand{\mcoeff}{\left| \kappa(\lambda) \right|}
\newcommand{\coeff}{\kappa(\lambda)}
\newcommand{\dive}{{\rm div}}
\newcommand{\I}[2]{I( #1, #2 )}

\newcommand{\di}{d_{\varepsilon, a}}

\newcommand{\qed}{$\Box$}
\newcommand{\diff}{\partial}
\newcommand{\ddu}{\tau_{u}}
\newcommand{\ddv}{\tau_{v}}
\newcommand{\ddz}{e_{z}}
\newcommand{\haus}{{\cal H}^{1}} 
\newcommand{\Om}{\Omega}
\makeatletter
\@addtoreset{equation}{section}
\makeatother

\newtheorem{theorem}{Theorem}[section]

\begin{document}
\null
\vspace{1cm}
\begin{center}
{\LARGE Local calibrations for minimizers of the  Mumford-Shah 

\vspace{.5cm}

 functional with rectilinear discontinuity sets}
\vspace{1cm}

 {\normalsize Gianni \sc Dal Maso}

\vspace{.25cm}
 {\normalsize Maria Giovanna \sc Mora}

\vspace{.25cm}
{\normalsize Massimiliano \sc Morini}

\vspace{1cm}

{\normalsize S.I.S.S.A.}\\
{\normalsize Via Beirut 2-4, 34014 Trieste, Italy}\\
{\normalsize e-mail: \tt dalmaso@sissa.it, mora@sissa.it, 
morini@sissa.it}

\vspace{3cm}

\begin{minipage}[t]{11.5cm}
{\normalsize
\centerline{ {\bf Abstract}}
\noindent Using a calibration method, 
we prove that, if $w$ is a function which satisfies
all Euler conditions for the Mumford-Shah functional on a two-dimensional
open set $\Omega$, and the discontinuity set of $w$ is a segment 
connecting two boundary points, then for every point $(x_0, y_0)$ of
$\Om$ there exists a neighbourhood $U$ of $(x_0, y_0)$ such that
$w$ is a minimizer of the Mumford-Shah functional on $U$ with respect 
to its own boundary values on $\partial U$.}
\vspace{.5cm}\\


\vspace{.5cm}
\end{minipage}
\end{center}

\setcounter{page}{0}
\thispagestyle{empty}

\vfill


\pagebreak

\clearpage
\section{Introduction}

This paper deals with Dirichlet minimizers of the Mumford-Shah 
functional (see \cite{Mum-Sha1} and \cite{Mum-Sha2})
\begin{equation} \label{g1}
\int_{\Om} |\nabla u(x,y)|^{2} dx\,dy +  \haus(S_{u}) \,,
\end{equation}
where $\Omega$ is a bounded open subset of $\rdue$ with a Lipschitz 
boundary, $\haus$ is the one-dimensional Hausdorff measure,
$S_{u}$ is the set of essential 
discontinuity points of the unknown function $u$, while $\nabla u$ 
denotes its approximate gradient (see \cite{Amb} or \cite{Amb-Fus-Pal}).

A {\it Dirichlet minimizer\/} of (\ref{g1}) in $\Om$  is a function $w$ which 
belongs to the space $SBV(\Om)$ of special functions of bounded 
variation in $\Om$ (see 
\cite{Amb} or \cite{Amb-Fus-Pal}) and satisfies the inequality
$$
\int_{\Om} |\nabla w(x,y)|^{2} dx \,dy + \haus(S_{w}) \,\le\,
\int_{\Om} |\nabla u(x,y)|^{2} dx \,dy + \haus(S_{u})
$$
for every function $u\in SBV(\Om)$ with the same trace as $w$ on $\partial\Om$. 

Suppose that $w$ is a Dirichlet minimizer 
of (\ref{g1}) in $\Om$ and that $S_{w}$ is a regular curve. Then the following 
equilibrium conditions are satisfied (see \cite{Mum-Sha1} and
\cite{Mum-Sha2}):
\begin{description}
\item[(a)] $w$ is harmonic on $\Om\setminus S_{w}$;
\item[(b)] the normal derivative of $w$ vanishes on both sides of $S_{w}$;
\item[(c)] the curvature of $S_{w}$ is equal to the
difference of the squares of the tangential derivatives of $w$
on both sides of $S_{w}$.
\end{description}
Elementary examples show that conditions (a), (b), and (c) are not 
sufficient for the Dirichlet minimality of $w$.

In this paper we prove that, if 
$S_{w}$ is a straight line segment connecting two 
points of $\partial \Om$, and  the tangential derivatives
$\diff_{\tau} w$ and $\diff_{\tau}^2 w$ 
of $w$ do not vanish on both sides of $S_{w}$, then (a), (b), and (c) 
imply that every point $(x_{0}, y_{0})$ in $\Om$ has
an open neighbourhood $U$ 
such that $w$ is a Dirichlet minimizer of (\ref{g1}) in $U$. In 
other words, under our assumptions, conditions (a), (b), and (c) are 
also sufficient for the Dirichlet minimality in small domains.
We hope that our proof will be useful in the future to achieve the same 
result without our special assumptions on $S_{w}$.

The proof is obtained by using the calibration method 
adapted in
\cite{Alb-Bou-DM}
to the functional (\ref{g1}). 
We construct an explicit calibration for $w$ 
in the cylinder $U{ \times }\erre$, where $U$ is a suitable 
neighbourhood of $(x_{0}, y_{0})$. This construction is elementary 
when $(x_{0}, y_{0}) \notin S_{w}$ (see \cite{Alb-Bou-DM}), so we 
consider only the case $(x_{0}, y_{0})\in S_{w}$.

The plan of the paper is the following. 
In Section 2 we fix the notations and we
recall the main result of \cite{Alb-Bou-DM}. 
In Theorem \ref{x-x} we consider
the special case of the function 
$$w(x,y) := \begin{cases}
\;\;\; x & \text{if $y>0$}, \\
-x & \text{if $y<0$},
\end{cases}$$
and give in full detail the first 
example of a calibration for a 
discontinuous function
which is not locally constant.
In Theorem \ref{x-x+1} we adapt the same construction to the function
$$w(x,y) := \begin{cases}
x+ 1 & \text{if $y>0$}, \\
x & \text{if $y<0$}.
\end{cases}$$
In Section 4
we consider the general case of a function $w$ satisfying (a), 
(b), and (c) and with $S_{w}=\{ (x,y) \in \Omega : y=0 \}$.
If $S_w$ is connected, only two situations are possible: 
\begin{eqnarray}
& \partial_{x} w(x,0+) = - \partial_{x} w(x,0-) \quad  \hbox{on} \ \, S_{w}, & 
\label{opp} \\
& \partial_{x} w(x,0+) = \partial_{x} w(x,0-) \quad  \hbox{on} \ \, S_{w}. & 
\label{conc}
\end{eqnarray}
The former case (\ref{opp}) is studied in Theorem \ref{harm} by a 
suitable change of variables and by adding two new parameters to the 
construction used in Theorem \ref{x-x}.
The minor changes for the case (\ref{conc})
are considered in Theorem \ref{harm2}.

\section{Preliminary results}

Let $\Omega$ be a bounded open subset of $\rdue$ with a Lipschitz 
boundary and let 
$$\Omega_{0} = \{ (x,y) \in \Omega : y \neq 0 \}, \qquad
S = \{ (x,y)\in \Omega : y=0 \}.$$
For every vector field $\varphi : 
\Om { \times } \erre \to \rdue { \times } \erre$ we define the maps
$\varphi^{x}, \; \varphi^{y}, \; \varphi^{z} : \Om { \times } \erre \to 
\erre$ by 
$$\varphi (x,y,z) = (\varphi^{x}(x,y,z), \varphi^{y}(x,y,z), 
\varphi^{z}(x,y,z)).$$
We shall consider the collection ${\cal F}$ of all piecewise $C^{0}$ 
vector fields $\varphi : 
\Om { \times } \erre \to \rdue { \times } \erre$ 
with the following property: there exists a finite 
number $g_{1}, \ldots, g_{k}$ of functions in $C^{1}(\overline{\Om})$ 
such that the sets
$$A_{i} : = \{ (x,y,z): (x,y)\in \Om, \ g_{i}(x,y) < z < g_{i+1}(x,y) \}$$
are nonempty and $\varphi \in C^{1}(\overline{A_{i}}, \rdue { 
\times } \erre)$ for $i=0, \ldots, k$, where we put $g_{0}=-\infty$ 
and $g_{k+1} = + \infty$.
Therefore, the discontinuity set of a vector field in ${\cal F}$ is contained in 
a finite number of regular surfaces.

Let $w \in C^{1}(\Om_{0})$ be a function such that
$\int_{\Om_{0}} |\nabla w|^{2} dx\, dy < + \infty$.
The upper trace of $w$ on $S$ is denoted by $w(x,0+)$, and the lower 
trace by $w(x,0-)$. Therefore, the approximate upper and lower limits 
$w^{+}(x,0)$ and $w^{-}(x,0)$
are given by 
$$w^{+}(x,0) = \max \{ w(x,0+), w(x,0-) \} \qquad \hbox{and}
\qquad w^{-}(x,0) = \min \{ w(x,0+), w(x,0-) \}.$$
A {\it calibration \/} for $w$ is a bounded vector field $\varphi 
\in {\cal F}$ which is continuous on the graph of $w$ and
satisfies the following properties:
\begin{description}
\item[(a)] $\dive \varphi = 0$ in the sense of distributions in 
$\Om { \times } \erre$;
\item[(b)] $(\varphi^{x}(x,y,z))^{2} + (\varphi^{y}(x,y,z))^{2} \leq 
4 \varphi^{z}(x,y,z)$ at every continuity point $(x,y,z)$ of $\varphi$;
\item[(c)]  $(\varphi^{x}, \varphi^{y})(x,y,w(x,y)) = 2 \nabla  w(x,y)$
and
$\varphi^{z}(x,y,w(x,y)) = |\nabla w(x,y)|^{2}$ for every $(x,y)\in 
\Om_{0}$;
\item[(d)] $\displaystyle \left( \int_{t_{1}}^{t_{2}} \varphi^{x}(x,y,z)\, dz 
\right)^{2} + \left( \int_{t_{1}}^{t_{2}} \varphi^{y}(x,y,z)\, dz 
\right)^{2} \leq 1$ for every $(x,y)\in \Om$ and for every $t_{1}, 
t_{2}\in \erre$;
\item[(e)] $\displaystyle \int_{w^{-}(x,0)}^{w^{+}(x,0)} \varphi^{x}(x,0,z)\, dz =0$ 
and
$\displaystyle \int_{w^{-}(x,0)}^{w^{+}(x,0)} \varphi^{y}(x,0,z)\, dz =1$ for every 
$(x,0)\in S$.
\end{description}

The following theorem is proved in \cite{Alb-Bou-DM}.

\begin{theorem}
If there exists a calibration $\varphi$ for $w$, then $w$ is a 
Dirichlet minimizer of the Mumford-Shah functional (\ref{g1}) in $\Om$.
\end{theorem}

If $\Om$ is a circle with centre on the $x$-axis, and
$w\in C^{1}(\Om_{0})$ with $\int_{\Om_{0}} |\nabla w|^{2} dx\, dy < + \infty$,
then $w$ satisfies the Euler conditions (a), (b), and (c)
if and only if $w$ has one of the following forms:
\begin{equation}\label{caso1}
w(x,y) =
\begin{cases}
\;\;\; u(x,y) & \text{if $y>0$}, \\
-u(x, y) + c_1 & \text{if $y<0$},
\end{cases}
\end{equation}
or
\begin{equation}\label{caso2}
w(x,y) =
\begin{cases}
u(x,y) + c_2 & \text{if $y>0$}, \\
u(x, y)  & \text{if $y<0$},
\end{cases}
\end{equation}
where $u\in C^1(\Om)$ is harmonic with normal derivative 
vanishing on $S$ and $c_1$, $c_2$ are real constants. 
For our purposes, it is enough to consider the case $c_1=0$ in (\ref{caso1})
and $c_2 =1$ in (\ref{caso2}).

\section{A model case}

In this section we consider in (\ref{caso1}) and in 
(\ref{caso2}) the particular function $u(x,y)=x$ and we
deal with the minimality of the functions
\begin{equation}\label{ics}
w (x,y) :=
\begin{cases}
\;\;\; x & \text{if $y>0$}, \\
-x & \text{if $y<0$},
\end{cases}
\end{equation}
and
\begin{equation}\label{costante}
w (x,y) :=
\begin{cases}
x+ 1 & \text{if $y>0$}, \\
x & \text{if $y<0$}.
\end{cases}
\end{equation}
The aim of the study of these simpler cases (but we will see
that they involve the main difficulties) is to clarify the ideas of the general construction.

\begin{theorem}\label{x-x}
Let $w: \rdue \to \erre$ be the 
function defined by
$$w(x,y) := \begin{cases} \;\;\;  x & \text{if $y > 0$}, \\
-x & \text{if $y < 0$}. 
\end{cases}$$
Then every point $(x_{0}, y_{0}) \neq (0,0)$ has an open 
neighbourhood $U$ such that $w$ is a Dirichlet minimizer 
in $U$ of
the Mumford-Shah functional (\ref{g1}).
\end{theorem}

{\bf Proof.}
The result follows by Theorem 4.1 of \cite{Alb-Bou-DM} if $y_{0}\neq 0$.
We consider now the case $y_{0}=0$, assuming for simplicity 
that $x_0 >0$. 
We will construct a local calibration of $w$ near $(x_{0}, 0)$.
Let us fix  $\varepsilon >0$ such that 
\begin{equation}\label{130}
0 < \varepsilon < \frac{x_{0}}{10}, \qquad 0< \varepsilon < 
\frac{1}{32}.
\end{equation}
For $0<\delta<\varepsilon$ we
consider the open rectangle
$$U := \{ (x,y)\in \rdue: |x-x_0| <\varepsilon, |y| < \delta \}$$
and the following subsets of $U{ \times } \erre$ (see Fig.~1)
\begin{figure}[p]
\begin{center} 
  \includegraphics[height=0.9\textheight]{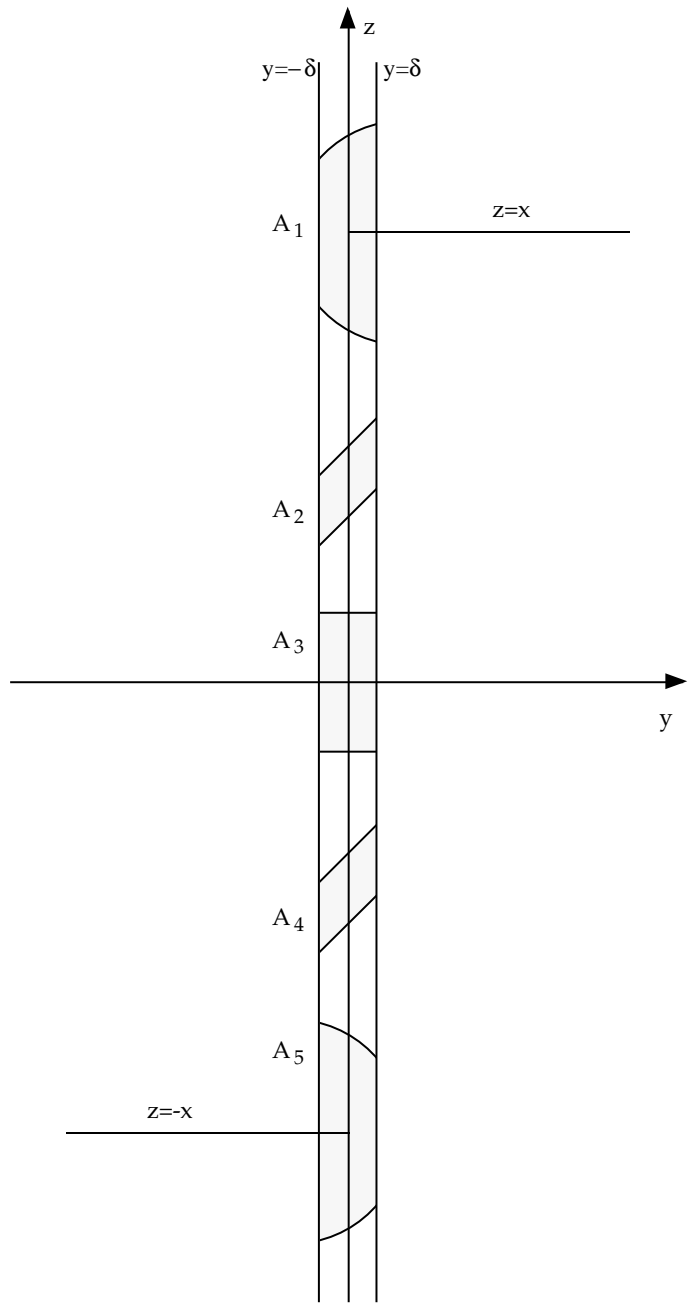}
\caption{Section of the sets 
$A_{1},\ldots,A_{5}$ at $x={\rm constant}$.}
\end{center}
\end{figure}
\begin{eqnarray*}
A_1 & := & \{(x,y,z)\in U{ \times } \erre : x-\alpha(y) < z < x+\alpha(y) \}, \\
A_2 & := & \left\{(x,y,z)\in U{ \times } \erre : b +\coeff \, y< z < b + \coeff \, y
+ h \right\}, \\ 
A_3 & := & \{(x,y,z)\in U{ \times } \erre : -h < z< h \}, \\
A_4 & := & \left\{(x,y,z)\in U{ \times } \erre : -b + \coeff \, y -h < z < -b
+ \coeff \, y \right\}, \\
A_5 & := & \{(x,y,z)\in U{ \times } \erre : -x-\alpha(-y) < z < 
-x+\alpha(-y) \}, \end{eqnarray*}
where 
$$\alpha(y) :=  \sqrt{4\varepsilon^2 - (\varepsilon-y)^2},$$
$$h:=\frac{x_0 -3\varepsilon}{4}, \qquad \coeff := 
\frac{\lambda}{4} - \frac{1}{\lambda},  \qquad
b := 2h+ \coeff \, \delta, \qquad \lambda := 
\frac{1-4\varepsilon}{2h}.$$
We will assume that
\begin{equation}\label{120}
\delta < \frac{x_0 -3\varepsilon}{8 \mcoeff},
\end{equation}
so that the sets $A_1, \dots, A_5$ are pairwise disjoint.

For every $(x,y,z) \in U{ \times } \erre $, let us 
define the vector 
$\varphi(x,y,z) = (\varphi^x, \varphi^y, \varphi^z)(x,y,z)\in \rtre$ as 
follows:
$$\begin{cases} 
\displaystyle \left(\frac{2(\varepsilon-y)}{\sqrt{(\varepsilon-y)^2 + 
(z-x)^2}}, \frac{-2(z-x)}{\sqrt{(\varepsilon-y)^2 + (z-x)^2}}, 1
\right) & \text{if $(x,y,z)\in A_1$}, \\
\\
\displaystyle \left( 0, \lambda, \frac{\lambda^2}{4} \right) & \text{if $(x,y,z)\in 
A_2$}, \\
\\
\displaystyle (f(y),0,1) & \text{if $(x,y,z)\in A_3$}, \\
\\
\displaystyle \left( 0, \lambda, \frac{\lambda^2}{4} \right) & \text{if $(x,y,z)\in 
A_4$}, \\ 
\\
\displaystyle \left( \frac{-2(\varepsilon+y)}{\sqrt{(\varepsilon+y)^2 + (z+x)^2}}, 
\frac{2(z+x)}{\sqrt{(\varepsilon+y)^2
+ (z+x)^2}}, 1 \right) & \text{if $(x,y,z)\in A_5$}, \\
\\
\displaystyle (0,0,1) & \text{otherwise}, 
\end{cases}$$
where 
$$f(y) := - \frac{1}{h}
\left( \int_{0}^{\alpha(y)} \frac{\varepsilon-y}{\sqrt{t^2 + (\varepsilon-y)^2}}\, dt -
\int_0^{\alpha(-y)} \frac{\varepsilon+y}{\sqrt{t^2 + (\varepsilon+y)^2}}\, dt
\right).$$

Note that $A_1 \cup A_5$ is an open neighbourhood of $\graph  
(w)\cap ( U { \times } \erre )$.
The purpose of the definition of $\varphi$ in $A_1$ and $A_5$ 
(see Fig.~2)
\begin{figure}[ht]
\begin{center}
  \includegraphics[width=14cm]{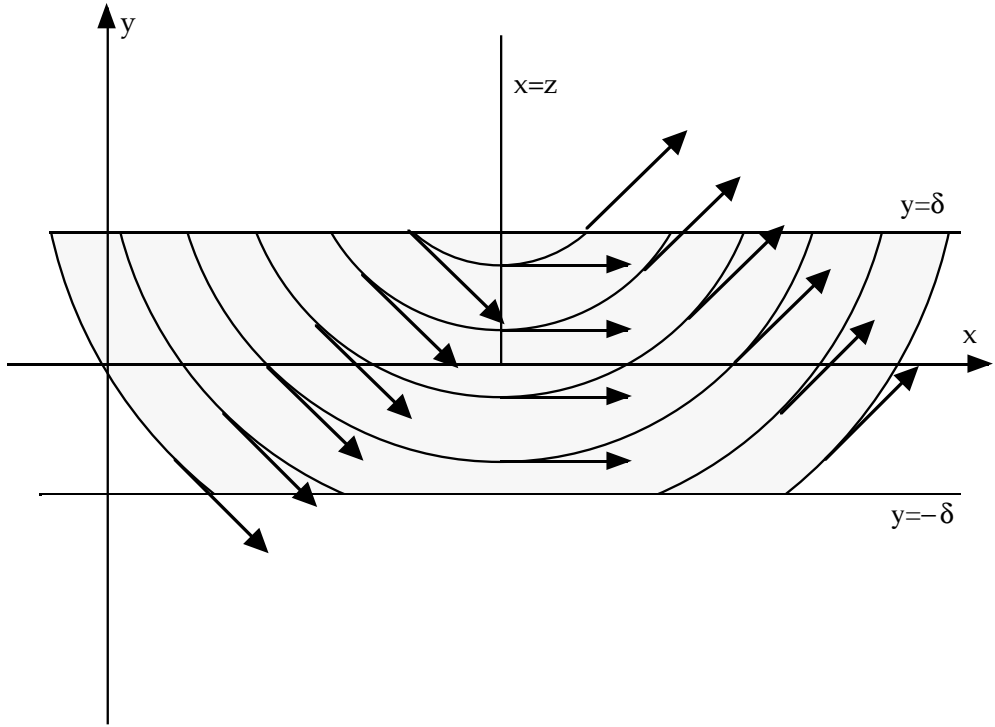}
\caption{Section of the set $A_{1}$ at $z={\rm constant}$.}
\end{center}
\end{figure}
is to provide a divergence free vector field satisfying 
condition (c) of Section~2 and such that
\begin{eqnarray*}
\varphi^y (x,0,z) > 0 & & \hbox{for} \;\; |z|<x, \\
\varphi^y (x,0,z) < 0 & & \hbox{for} \;\; |z|>x.
\end{eqnarray*}
These properties are crucial in order to obtain (d) and (e) simultaneously.

The r\^ole of $A_2$ and $A_4$ is to 
give the main contribution to the integral in (e). 
To explain this fact, suppose, for a moment, that $\varepsilon = 0$; 
in this case we would have $A_1= A_5 = \emptyset$ and
$$\int_{-x}^x \varphi^y (x,0,z)\, dz = 1,$$
so that the $y$-component of equality (e) would be satisfied.

The purpose of the definition of $\varphi$ in $A_3$ is to correct the 
$x$-component of $\varphi$, in order to obtain (d).

\

We shall prove that, for a suitable choice of 
$\delta$, the vector
field $\varphi$ is a calibration for $w$ in the rectangle $U$.

Note that for a given $z\in \erre$ we have
\begin{equation}\label{div2}
\partial_x \varphi^x (x,y,z) + \partial_y \varphi^y (x,y,z)=0
\end{equation}
for every $(x,y)$ such that $(x,y,z)\in A_1 \cup A_5$. This implies
$\varphi$ is divergence free in $A_1 \cup A_5$.
Moreover $\dive \varphi =0$ in the other sets $A_i$,
and the normal component  of 
$\varphi$  is continuous across $\partial A_i$:
the choice of $\coeff$ ensures that this property holds for $\partial 
A_{2}$ and $\partial A_{4}$ (see Fig.~3).
\begin{figure}[hp]
\begin{center} 
  \includegraphics[height=0.9\textheight]{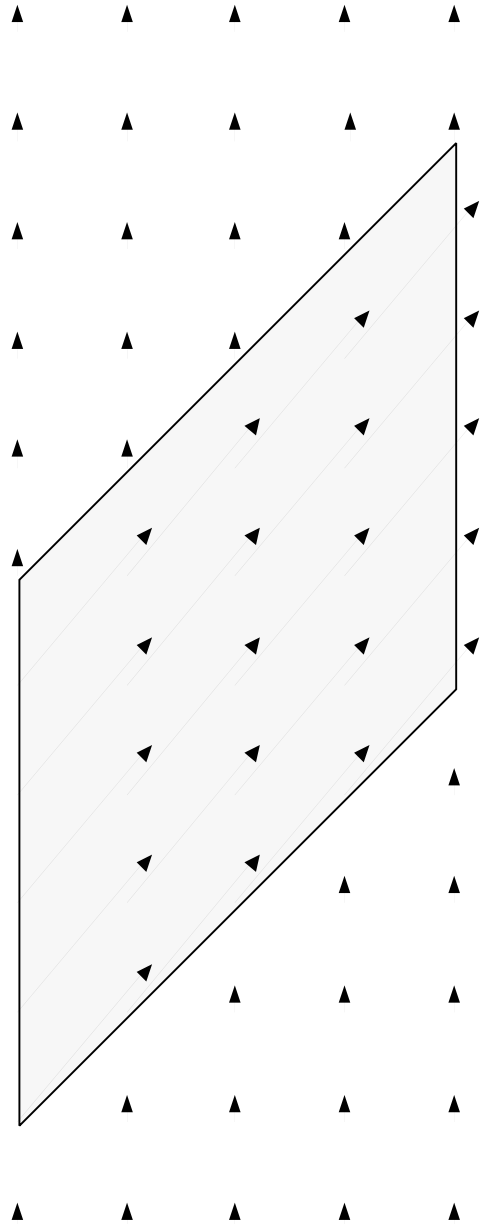}
  \caption{Section of the set $A_{2}$ at $x={\rm constant}$.}
\end{center}
\end{figure}
Therefore
$\varphi$ is divergence free in the sense of distributions
in $U{ \times } \erre$.

On the graph of $w$ we have
$$\varphi(x,y,w(x,y)) = \begin{cases}
(\;\;\; 2,0,1) & \text{if $y>0$}, \\
(-2,0,1) & \text{if $y<0$},
\end{cases}$$
so condition (c) is satisfied.

Inequality (b) is clearly satisfied in all regions: the only non 
trivial case is $A_3$, 
where we have, using (\ref{130}),
$$|f(y)|  \leq \frac{ 4\left( \alpha(y) + 
\alpha(-y)
\right) }{x_0 -3\varepsilon} 
 \leq  \frac{8\sqrt{3}\varepsilon }{x_0 - 3\varepsilon}
<2.$$

We now compute
\begin{equation}
\int_{-x}^{x} \varphi^y(x,y,z)\, dz .
\end{equation}
Let us fix $y$ with $|y|< \delta$. Since $\varphi^y(x,y,z)$ depends  
on $z-x$, we have
\begin{equation}\label{01}
\int_{x-\alpha(y)}^x \varphi^y(x,y,z)\, dz = \int_x^{x+\alpha(y)}
\varphi^y(\xi, y, x)\, d\xi.
\end{equation}
Using (\ref{div2}) and applying the divergence 
theorem to the curvilinear triangle 
$$T=\{ (\xi, \eta) \in \rdue: \xi> x, \; \eta< y, \; (\varepsilon-\eta)^2 + 
(x-\xi)^2<
4\varepsilon^{2} \}$$
(see Fig.~4), 
we obtain
\begin{equation}\label{02}
\int_x^{x+\alpha(y)}
\varphi^y(\xi, y, x)\, d\xi= \int_{-\varepsilon}^y 
\varphi^x(x, \eta, x)\, d\eta = 2(y+\varepsilon).
\end{equation}
\begin{figure}[t]
\begin{center}
  \includegraphics[width=14cm]{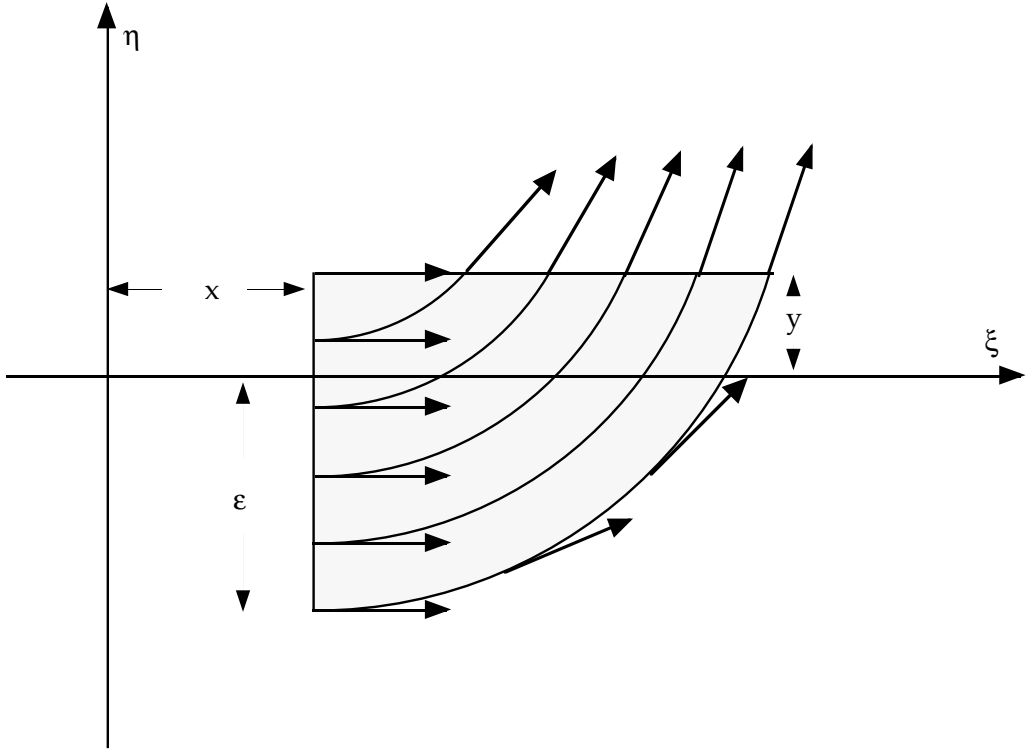}
  \caption{The curvilinear triangle $T$.}
\end{center}
\end{figure} 
{}From (\ref{01}) and (\ref{02}), we get
\begin{equation}\label{03}
\int_{x-\alpha(y)}^x \varphi^y(x,y,z)\, dz =
2(y+\varepsilon).
\end{equation}
Similarly we can prove that
\begin{equation}\label{04}
\int_{-x}^{-x +\alpha(-y)}\varphi^y(x,y,z)\, dz =
2(-y+\varepsilon).
\end{equation}
Using the definition of $\varphi$ in $A_2$, $A_3$, $A_4$, we obtain
\begin{equation}\label{phiy}
\int_{-x}^x \varphi^y(x,y,z)\, dz = 1.
\end{equation}
On the other hand, by the definition of $f$, we have immediately that
\begin{equation}\label{phix}
\int_{-x}^x \varphi^x(x,y,z)\, dz = 0.
\end{equation}
{}From these equalities it follows in particular 
that condition (e) is satisfied on the jump set 
$S_{w} \cap U = \{ (x,y)\in U: y=0\}$.

Let us begin now the proof of (d).
Let us fix $(x,y)\in U$.
For every $t_{1}< t_{2}$ we set
$$\I{t_{1}}{t_{2}} := \int_{t_{1}}^{t_{2}} (\varphi^x,\varphi^y)(x,y,z)\, dz.$$
It is enough to consider the case $-x-\alpha(-y) \leq t_{1} \leq t_{2} 
\leq x - \alpha(y)$. We can write
\begin{eqnarray*}
\I{t_{1}}{t_{2}} & = & \I{t_{1}}{-x} + \I{-x}{x} + \I{x}{t_{2}}, \\
\I{t_{1}}{-x} & = & \I{t_{1} \land (-x +\alpha(-y))}{-x} + \I{t_{1} 
\lor (-x +\alpha(-y))}{-x +\alpha(-y)}, \\
\I{x}{t_{2}} & = & \I{x}{t_{2} \lor (x-\alpha(y))} + \I{x-\alpha(y)}{
t_{2} \land (x -\alpha(y))}.
\end{eqnarray*}
Therefore
$$\I{t_{1}}{t_{2}}  =  \I{-x}{x} + \I{t_{1} \land (-x 
+\alpha(-y))}{-x} + \I{x}{t_{2} \lor (x-\alpha(y))}$$
\begin{equation}\label{111}
\mbox{} + \I{t_{1} 
\lor (-x +\alpha(-y))}{t_{2} \land (x -\alpha(y))} - 
\I{-x+\alpha(-y)}{x-\alpha(y)}. 
\end{equation}

Let $B$ be the ball of radius $4\varepsilon$ 
centred at $(0,-4\varepsilon)$. We want to prove that
\begin{equation}\label{109}
\I{x}{t} \in \overline{B}
\end{equation}
for every
$t$ with $x-\alpha(y) \leq t \leq x+\alpha(y)$.
Let us denote the components of $\I{x}{t}$ by 
$a^x$ and  $a^y$.
Arguing as in the proof of (\ref{03}), we get the identity
\begin{equation}
a^{y} = 2(\varepsilon-y) - 2 \sqrt{(t-x)^2 + (\varepsilon-y)^2} \leq 0.
\end{equation}
As $|\varphi^x|\leq 2$, we have also
$$( a^{x} )^2 \leq 4(t-x)^2
= ( 2(\varepsilon -y) - a^{y}
)^2 - 4(\varepsilon -y)^2.$$
{}From these estimates it follows that
$$( a^{x} )^2 +
( a^{y} + 4\varepsilon
)^2 \leq 16 \varepsilon^2,$$
which proves (\ref{109}).
In the same way we can prove that
\begin{equation}\label{110}
\I{t}{-x} \in \overline{B}
\end{equation}
for every $t$ with $-x-\alpha(-y) \leq t \leq -x +\alpha(-y).$

If $f(y)\geq 0$, we define
$$C : = ([ 0, 2h f(y) ] { \times } [0, 
{\textstyle \frac{1}{2}}
-2\varepsilon]) \cup (\{ 2hf(y) \} {\times }  [0, 1-4\varepsilon]);$$
if $f(y)\leq 0$, we simply replace $[ 0, 2h f(y) ]$ by $[ 2h f(y), 
0 ]$.
From the definition of $\varphi$ in $A_{2}$, $A_{3}$, $A_{4}$, it 
follows that
\begin{equation}\label{112}
\I{-x+\alpha(-y)}{x-\alpha(y)} = (2h f(y), 1-4\varepsilon)
\end{equation}
and
\begin{equation}\label{114}
\I{s_{1}}{s_{2}} \in C
\end{equation}
for $-x+\alpha(-y) \leq s_{1} \leq s_{2} \leq x - \alpha(y)$.
Let $D := C - (2h f(y), 1-4\varepsilon)$, i.e.,
$$D = ([ -2h f(y) , 0 ] { \times } [-1 + 
4\varepsilon, 
-{\textstyle \frac{1}{2}}
+2\varepsilon]) \cup (\{ 0 \} {\times }  [-1+4\varepsilon , 0]),$$
for $f(y)\geq 0$; the interval $[ -2h f(y) , 0 ]$ is replaced by
$[ 0, -2h f(y) ]$ when $f(y)\leq 0$.
From (\ref{111}), (\ref{phiy}), (\ref{phix}),
(\ref{109}), (\ref{110}), (\ref{112}) and (\ref{114})
we obtain 
\begin{equation}\label{122}
\I{t_{1}}{t_{2}} \in (0,1) + 2\overline{B} + D.
\end{equation}
As $f(0)=0$, we can choose $\delta$ so that (\ref{120}) is satisfied 
and
\begin{equation}\label{121}
|2h f(y)| = \frac{x_{0} -3\varepsilon}{2} |f(y)| \leq \varepsilon
\end{equation}
for $| y | < \delta$.
It is then easy to see that, by (\ref{130}), 
the set $(0,1) + 2\overline{B} + D$ is contained in the unit ball centred at 
$(0,0)$.
So that (\ref{122}) implies (d). \qed

\

{\bf Remark.} The assumption $(x_0 , y_0)\neq (0,0)$ 
in Theorem \ref{x-x} cannot be dropped. 
Indeed, there is no 
neighbourhood $U$ of $(0,0)$ such that $w$ is a Dirichlet minimizer of 
the Mumford-Shah 
functional in $U$.

To see this fact, let $\psi$ be a function defined on the square
$Q = (- 1 ,1) { \times } (-1 ,1)$ satisfying the boundary condition
$\psi = w$ on $\partial Q$
and such that $S_{\psi} = ((-1, -1/2) \cup (1/2, 1)) { \times } \{ 0 \}$.
For every $\varepsilon$, let $\psi_{\varepsilon }$ be the function
defined on $Q_{\varepsilon} = \varepsilon Q$ by
$\psi_{\varepsilon}(x, y) := \varepsilon \psi (x / \varepsilon,
y / \varepsilon)$.
Note that $\psi_{\varepsilon}$ satisfies the boundary condition
$\psi_{\varepsilon}= w$ on $\partial Q_{\varepsilon}$.
Let us
compute the Mumford-Shah functional for $\psi_{\varepsilon}$ on $Q_{\varepsilon}$:
$$\int_{Q_{\varepsilon}} |\nabla \psi_{\varepsilon}|^2 dx\, dy + 
\haus(S_{\psi_{\varepsilon}}) =
\varepsilon^2 \int_{Q} |\nabla \psi |^2 dx\, dy + \varepsilon.$$
Since 
$$\int_{Q_{\varepsilon}} |\nabla w|^2 dx\, dy + \haus(S_{w}) = 4\varepsilon^2 + 
2\varepsilon,$$ 
we have 
$$\int_{Q_{\varepsilon}} |\nabla \psi_{\varepsilon}|^2 dx\, dy + 
\haus(S_{\psi_{\varepsilon}}) 
< \int_{Q_{\varepsilon}} |\nabla w|^2 dx\, dy + \haus(S_{w})$$
for $\varepsilon$ 
sufficiently small. 

\

The construction shown in the proof of Theorem \ref{x-x} can be easily adapted to 
define a calibration for the function $w$ in (\ref{costante}).

\begin{theorem}\label{x-x+1}
Let $w: \rdue \to \erre$ be the 
function defined by
$$w(x,y) : = \begin{cases}
x+1  & \text{if $y > 0$}, \\
x & \text{if $y < 0$}. 
\end{cases}$$
Then every point $(x_{0}, y_{0})\in \rdue$ has an open 
neighbourhood $U$ such that $w$ is a Dirichlet minimizer in $U$ of
the Mumford-Shah functional (\ref{g1}).
\end{theorem}

{\bf Proof.}
The result follows by Theorem 4.1 of \cite{Alb-Bou-DM} if $y_{0}\neq 0$.
We consider now the case $y_{0}=0$;
we will construct a local calibration of $w$ near $(x_{0}, 0)$,
using the same technique as in Theorem \ref{x-x}.
We give only the new definitions of the sets $A_1, \ldots, A_5$ and of the 
function $\varphi$, and leave to the reader the verification of the fact that
this function is a calibration for suitable values of the involved
parameters.

Let us fix  $\varepsilon >0$ such that 
\begin{equation}\label{130b}
0 < \varepsilon < \frac{1}{24}, \qquad 0< \varepsilon < 
\frac{1}{32}.
\end{equation}
For $0<\delta<\varepsilon$ we
consider the open rectangle
$$U := \{ (x,y)\in \rdue: |x-x_0| <\varepsilon, |y| < \delta \}$$
and the following subsets of $U { \times } \erre$
\begin{eqnarray*}
A_1 & := & \{(x,y,z)\in U{ \times } \erre  : x+1 - \alpha(y) < z < x+1 +\alpha(y) \}, \\
A_2 & := & \{(x,y,z)\in U {\times } \erre : b +  \kappa (\lambda ) \, y + 3h < z 
< b + \kappa (\lambda ) \, y + 4h  \}, \\ 
A_3 & := & \{(x,y,z)\in U { \times } \erre : x_0 + 3 \varepsilon + 2h <z < x_0 +
3 \varepsilon +3 h \}, \\
A_4 & := & \{ (x,y,z)\in U { \times } \erre : b + \kappa (\lambda ) \, y < z <
b + \kappa (\lambda ) \, y + h \}, \\
A_5 & := & \{(x,y,z)\in U { \times } \erre: x - \alpha (-y) <z < x+ \alpha (-y)   \}, \\
\end{eqnarray*}
where 
$$\alpha(y) :=  \sqrt{4\varepsilon^2 - (\varepsilon-y)^2},$$
$$h :=\frac{1 -6 \varepsilon}{5}, \qquad \kappa ( \lambda ) := 
\frac{\lambda}{4} - \frac{1}{\lambda},  \qquad
b : = x_0 + 3 \varepsilon + \kappa (\lambda ) \, \delta, \qquad 
\lambda := 
\frac{1- 4 \varepsilon }{2h}.$$
We will assume that
\begin{equation}\label{120b}
\delta < \frac{1 -6 \varepsilon}{10 | \kappa (\lambda ) | },
\end{equation}
so that the sets $A_1, \ldots, A_5$ are pairwise disjoint.

For every $(x,y,z) \in U { \times } \erre$, let us define the vector
$\varphi ( x,y,z) \in \rtre$ as follows:
$$\begin{cases} 
\displaystyle \left(\frac{2( \varepsilon -y)}{\sqrt{( \varepsilon -y)^2 + (z-x-1)^2}},
\frac{-2(z-x-1)}{\sqrt{(\varepsilon -y)^2 + (z-x-1)^2}}, 1
\right) & \text{if $(x,y,z)\in A_1$}, \\
\\
\displaystyle \left( 0, \lambda, \frac{\lambda^2}{4} \right) & \text{if $(x,y,z)\in A_2$}, \\
\\
\displaystyle (f(y),0,1) & \text{if $(x,y,z)\in A_3$}, \\
\\
\displaystyle \left( 0, \lambda, \frac{\lambda^2}{4} \right) & \text{if $(x,y,z)\in A_4$}, \\
\\
\displaystyle \left( \frac{2(\varepsilon +y)}{\sqrt{(\varepsilon +y)^2 + (z-x)^2}}, \frac{2(z-x)}{\sqrt{(
\varepsilon +y)^2 + (z-x)^2}}, 1 \right) &
\text{if $(x,y,z)\in A_5$}, \\
\\
\displaystyle (0,0,1) & \text{otherwise}, 
\end{cases}$$
where
$$f(y) := - \frac{2}{h}\left( \int_{0}^{
\alpha(y)} \frac{\varepsilon -y}{\sqrt{t^2 + (\varepsilon -y)^2}} \, dt
+ \int_{0}^{
\alpha(-y)} \frac{\varepsilon +y}{\sqrt{t^2 + (\varepsilon +y)^2}}\, dt
\right)$$
for every $|y| <  \delta$.
\qed

\section{The general case}

In this section we denote by $\Omega$ a ball in $\rdue$ centred at $(0,0)$
and we consider as $u$ in (\ref{caso1}) and in (\ref{caso2}) a generic 
harmonic function with normal derivative vanishing on $S$.
We add the technical assumption that the 
first and second order tangential derivatives of 
$u$ are not zero on $S$.

\begin{theorem}\label{harm}
Let $u: \Omega \to \erre$ be a harmonic function such that
$\partial_{y}u(x,0)=0$ for $(x,0)\in\Omega$, and let 
$w: \Omega \to \erre$ be the function defined by
$$
w(x,y) := \begin{cases}\;\;\;u(x,y) & \text{for $y>0$},
\\
-u(x, y) & \text{for $y<0$}.
\end{cases}$$
Assume that $u_{0}:=u(0,0)\neq 0$, $\partial_{x}u(0,0)\neq 0$, and
$\partial_{x}^2 u(0,0)\neq 0$. 
Then there 
exists an open 
neighbourhood $U$ of $(0,0)$ such that $w$ is a Dirichlet minimizer in $U$ of
the Mumford-Shah functional (\ref{g1}).
\end{theorem}

{\bf Proof.} We may assume $u(0, 0)>0$ and $\partial_{x}u(0,0) >0$.
We shall give the proof only for $\partial_{x}^2 u(0,0)> 0$,
and we shall explain at the end the modification needed for
$\partial_{x}^2u(0,0)< 0$.
Let $v:\Omega\to \erre$ be the harmonic conjugate of $u$ that 
vanishes on $y=0$, i.e., the function satisfying 
$\partial_{x}v(x,y) = -\partial_{y}u(x,y)$, 
$\partial_{y}v (x,y) = \partial_{x}u(x,y)$, and
$v(x,0)=0$.

Consider a small neighbourhood $U$ of $(0,0)$ such that the map
$\Phi(x,y):= (u(x,y), v(x,y))$ is invertible on $U$ and $\diff_x u 
>0$ on $U$. We call
$\Psi$ the inverse function $(u,v) \mapsto (\xi(u,v), \eta(u,v))$, 
which is defined 
in the 
neighbourhood $V:= \Phi (U)$ of $(u_{0},0)$.
Note that, if $U$ is small enough, then $\eta (u,v) =0$  
if and only if $v=0$. Moreover, 
\begin{equation}\label{psijacob}
D\Psi = \left( \begin{array}{cc}
\diff_u \xi & \diff_v \xi  \\
\diff_u \eta  & \diff_v \eta  
\end{array} \right) =
\frac{1}{|\nabla u|^2} \left(
\begin{array}{cc}
\diff_x u & \diff_x v  \\
\diff_y u  & \diff_y v
\end{array}
\right),
\end{equation}
where, in the last formula, all functions are computed at $(x,y)= \Psi(u,v)$,
and so
$\diff_u \xi = \diff_v \eta$, $\diff_v \xi = - \diff_u \eta$ and
$\diff_u \eta (u,0)=0$, $\diff_v \eta (u,0) > 0$. In particular, $\xi$ 
and $\eta$ are harmonic, and
\begin{equation}\label{nulla}
\diff_{u}^{2}\eta (u,0) =0, \qquad \diff_{v}^{2}\eta (u,0) =0.
\end{equation}
On $U$ we will use the coordinate system $(u,v)$ given by $\Phi$.
By (\ref{psijacob}) the canonical 
basis of the tangent space to $U$ at a point $(x,y)$
is given by
\begin{equation}\label{diff}
\ddu = \frac{\nabla u}{|\nabla u|^2}, \qquad \ddv =\frac{\nabla v}{|\nabla v|^2}.
\end{equation}
For every $(u,v)\in V$,
let $G(u,v)$ be the matrix associated with the first fundamental form
of $U$ in the coordinate system $(u,v)$,
and let $g(u,v)$  be its determinant. 
By (\ref{psijacob}) and (\ref{diff}), 
\begin{equation}\label{gi}
g = ((\diff_u \eta)^2 + (\diff_v \eta)^2)^{2} = \frac{1}{|\nabla u 
(\Psi) |^4}.
\end{equation}
We set $\gamma(u,v) := \sqrt[4]{g(u,v)}$.

The calibration $\varphi(x,y,z)$ on $U { \times } \erre$ 
will be written as 
\begin{equation}\label{calibra2}
\varphi(x,y,z) = \frac{1}{\gamma^2 (u(x,y), 
v(x,y))}\phi(u(x,y), 
v(x,y), z).
\end{equation}
We will adopt the following representation for $\phi : V { \times } \erre \to \rtre$:
\begin{equation}\label{calibra}
\phi (u,v,z)=
\phi^u (u,v,z) \ddu +  \phi^v (u,v,z) \ddv +  \phi^z (u,v,z) 
\ddz,
\end{equation}
where $\ddz$ is the third vector of the canonical basis of $\rtre$, 
and $\ddu$, $\ddv$ are computed at the point $\Psi (u,v)$.
We now reformulate the conditions of Section 2
in this new coordinate system.
It is known from Differential Geometry (see, e.g., \cite[Proposition 
3.5]{Cha}) that,
if $X = X^u \ddu + X^v \ddv$ 
is a vector field on $U$,
then the divergence of $X$ is given by
\begin{equation}\label{div}
\dive X = \frac{1}{\gamma^2} (\diff_{u}(\gamma^2 X^u) +
\diff_{v} (\gamma^2 X^v)).
\end{equation}
Using (\ref{diff}), (\ref{gi}), (\ref{calibra2}), (\ref{calibra}), and
(\ref{div})
it turns out
that $\varphi$ is a calibration if
the following conditions are satisfied: 
\begin{description}
\item[(a)] $\displaystyle 
     \diff_{u} \phi^u +
	\diff_{v} \phi^v +
	\diff_{z} \phi^z = 0$ for every $(u,v,z) \in V { \times } \erre$;
\item[(b)] $\displaystyle (\phi^u (u,v,z))^2 + (\phi^v(u,v,z))^2 \leq 
            4 \phi^z(u,v,z)$
	for every  $(u,v,z) \in V { \times } \erre$;
\item[(c)] $\displaystyle \phi^u (u,v, \pm u)=\pm 2$,
	$\displaystyle \phi^v(u,v, \pm u)=0$, and 
	$\displaystyle \phi^z(u,v, \pm u)= 1$ 
	for every $(u,v)\in V$;
\item[(d)] $\displaystyle \left( \int_{t_1}^{t_2} \phi^u(u,v,z)\,
	dz \right)^2 + \left( \int_{t_1}^{t_2}
	\phi^v(u,v,z)\, dz \right)^2 \leq
	\gamma^2(u,v)$
	for every $(u,v) \in V$, $t_1, t_2 \in \erre$;
\item[(e)] $\displaystyle \int_{-u}^{u} \phi^u(u,0,z)\, dz = 0$ and 
	$\displaystyle \int_{-u}^u \phi^v (u,0,z)\, dz = \gamma(u,0)$
	for every $(u,0) \in V$.
\end{description}

Given suitable parameters $\varepsilon >0$, $h>0$, $\lambda >0$, that 
will be chosen later,
and assuming 
\begin{equation}\label{intorno}
V=
\{ (u,v): |u-u_{0}| < \delta, |v| < \delta \},
\end{equation}
with $\delta < \varepsilon$,
we consider the following subsets of $V{ \times } \erre$
\begin{eqnarray*}
A_1 & := & \{(u,v,z)\in V{ \times } \erre : u
-\alpha(v) < z < u +\alpha(v) \}, \\
A_2 & := & \left\{(u,v,z)\in V{ \times } \erre : 3h +\beta (u,v) < z < 
3h + 
\beta (u,v)
+ 1/\lambda \right\}, \\ 
A_3 & := & \{(u,v,z)\in V{ \times } \erre : -h < z< h \}, \\
A_4 & := & \left\{(u,v,z)\in V{ \times } \erre : -3h + \beta (u,v) 
-1/\lambda < z < -3h
+ \beta (u,v) \right\}, \\
A_5 & := & \{(u,v,z)\in V{ \times } \erre : -u -\alpha(-v) < z < -
u + \alpha(-v) \},
\end{eqnarray*}
where
$$\alpha(v) :=  \sqrt{4\varepsilon^2 - (\varepsilon -v)^2},$$
and $\beta$ is a suitable smooth function satisfying
$\beta (u,0) =0$, which will be defined later.
It is easy to see that, if $\varepsilon$ and $h$ are sufficiently 
small, while $\lambda$ is sufficiently large, then 
the sets $A_{1}, \ldots , A_{5}$ are pairwise disjoint,
provided $\delta$ is small enough. Moreover, 
since $\gamma(u,0)=
\diff_{v} \eta (u,0)>0$,
by continuity we may assume that
\begin{equation}\label{gamma}
\gamma(u,v) > 128\varepsilon \qquad \hbox{and}
\qquad \diff_{v}\eta (u,v) > 8\varepsilon 
\end{equation}
for every $(u,v)\in V$.

For $(u,v) \in V$ and $z\in\erre$ the vector 
$\phi(u,v,z)$
introduced in (\ref{calibra2}) is defined
as follows:
$$
\begin{cases}
\displaystyle \frac{2(\varepsilon- v)}{\sqrt{(\varepsilon-v)^2+(z-u)^2}} \ddu
-  \frac{2(z-u)}{\sqrt{(\varepsilon-v)^2
+(z-u)^2}} \ddv + \ddz
& \text{in $A_1$}, \\
\\
\displaystyle -\lambda \sigma (u,v) \frac{v}{\sqrt{(u-a)^2 + v^2}} \ddu +
\lambda \sigma (u,v) \frac{u-a}{\sqrt{(u-a)^2 + v^2}} \ddv +
\mu \ddz  
& \text{in $A_2$}, \\
\\
\displaystyle  f(v)\ddu + \ddz  & \text{in 
$A_3$}, \\
\\
\displaystyle -\lambda \sigma (u,v) \frac{v}{\sqrt{(u-a)^2 + v^2}} \ddu +
\lambda  \sigma (u,v) \frac{u-a}{\sqrt{(u-a)^2 + v^2}} \ddv +
\mu \ddz  
& \text{in $A_4$}, \\
\\
\displaystyle - \frac{2(\varepsilon+ v)}{\sqrt{(\varepsilon+ v)^2 + (z+ u)^2}} \ddu
+ \frac{2(z+ u)}{\sqrt{
(\varepsilon+ v)^2 + (z+ u)^2}} \ddv + \ddz  & 
\text{in $A_5$}, \\
\\
\displaystyle \ddz  & \text{\hphantom{in $A_{5}$,}\llap{otherwise,}}
\end{cases}$$
where 
\begin{eqnarray}\label{param}
& a < u_{0}- 11 \delta, \qquad \mu >0 & \\
\nonumber \\
& \displaystyle f(v) := - \frac{1}{h} \left( \int_0^{\alpha(v) }
\frac{(\varepsilon - v)}{\sqrt{t^2 + (\varepsilon - v)^2}}\, dt - 
\int_0^{\alpha(-v) }
\frac{(\varepsilon + v)}{\sqrt{t^2 + (\varepsilon + v)^2}}\, dt 
\right), & \nonumber \\
\nonumber \\
\label{sigma} \displaystyle
 & \sigma (u,v) := \frac{1}{2}
\gamma ( a + \sqrt{(u-a)^2 + v^2},0) 
- 2 \varepsilon. & 
\end{eqnarray}
We choose $\beta$ as the solution of the 
Cauchy problem
\begin{equation}\label{cauchy}
\begin{cases}
\lambda \sigma(u,v) (- v\, \partial_{u} \beta + (u-a) 
\partial_{v} \beta ) 
= (\mu -1) \sqrt{(u-a)^{2} + v^{2}}, & \\
\\
\beta (u,0)=0. & 
\end{cases}
\end{equation}
Since the line $v=0$ is not characteristic for the equation near 
$(u_{0},0)$,
there exists a unique solution $\beta \in C^{\infty}(V)$, provided 
$V$ is small enough.

In the coordinate system $(u,v)$ the definition of the field $\phi$
in $A_{1}$, $A_{3}$, and $A_{5}$
is the same as the definition of $\varphi$ in the proof of Theorem 
\ref{x-x}. The crucial difference is in the definition on the sets 
$A_{2}$ and $A_{4}$, where now we are forced to introduce two new parameters $a$ 
and $\mu$. Note that the definition given in Theorem \ref{x-x} can be
regarded as the limiting case as $a$ tends to $+\infty$.

By direct computations it is easy to see that $\phi$ satisfies condition 
(a) on $A_{1}$ and $A_{5}$.
Similarly,
the vector field
$$\left( -\frac{v}{\sqrt{(u-a)^2 + v^2}},
\frac{u-a}{\sqrt{(u-a)^2 + v^2}} \right)$$
is divergence free; 
since 
$(u-a)^{2}+v^{2}$ is constant along the integral curves of this field,
by construction 
the same property holds for $\sigma$, so that
$\phi$ satisfies condition (a) in $A_{2}$ and $A_{4}$.

In $A_3$, condition (a) is trivially satisfied.

Note that
the normal component
of $\phi$  is continuous across each $\partial A_i$: 
for the region $A_3$ this continuity is guaranteed by our
choice of $\beta$. This implies that (a) is satisfied in the sense of 
distributions on $V { \times } \erre$.

In order to satisfy condition (b),
it is enough to take
the parameter $\mu$ such that
$$\frac{\lambda^2}{4} 
\sigma^2(u,v)
\leq \mu$$
for every $(u,v)\in V$,
and require that
\begin{equation}\label{effe}
|f(v)|\leq 2.
\end{equation}
Since
\begin{equation}\label{stima4}
|f(v)|  \leq  \frac{\alpha(v) + \alpha(-v)}{h} 
\leq \frac{4\, \varepsilon}{h},
\end{equation}
inequality (\ref{effe}) is true if we impose
$$2 \, \varepsilon \leq h.$$

Looking at the definition of $\phi$ on $A_{1}$ and $A_{5}$, 
one can check that
condition (c) is satisfied.

Arguing as in the proof of (\ref{03}),
(\ref{04}), (\ref{phix}) in Theorem \ref{x-x}, 
we find that for every $(u,v) \in V$
$$\int_{-u}^{-u+\alpha(-v)} \phi^u(u,v,z)\, dz
+ \int_{-h}^h \phi^u(u,v,z)\, dz
+ \int_{u-\alpha(v)}^u \phi^u(u,v,z)\, dz =0,$$
$$\int_{-u}^{-u+\alpha(-v)} \phi^v(u,v,z)\, dz
+ \int_{-h}^h \phi^v(u,v,z)\, dz
+ \int_{u-\alpha(v)}^u \phi^v(u,v,z)\, dz = 4\varepsilon.$$
Now, it is easy to see that
\begin{equation}\label{phiu}
\int_{-u}^u \phi^u(u,v,z)\, dz = -2\sigma(u,v)
\frac{v}{\sqrt{(u-a)^2+v^2}},
\end{equation}
\begin{equation}\label{phiv}
\int_{-u}^u \phi^v(u,v,z)\, dz = 4\varepsilon
+ 2 \sigma(u,v)
\frac{u-a}{\sqrt{(u-a)^2+v^2}};
\end{equation}
since for $v=0$ we have
$$\sigma (u,0) = \frac{1}{2} \gamma(u,0)- 2\varepsilon,$$
condition (e) is satisfied.

By continuity, if $\delta$ is small enough,
we have
\begin{equation}\label{stima2}
\int_{-u}^u \phi^v(u,v,z)\, dz > \frac{7}{8} \gamma(u,v)
\end{equation}
for every $(u,v)\in V$.

From now on, we regard the pair $(\phi^u, \phi^v)$ as a vector in $\rdue$.
To prove condition (d) we set
$$I_{\varepsilon,a}(u,v,s,t) :=\int_s^t (\phi^u,\phi^v)(u,v,z)\, dz$$
for every $(u,v)\in V$, and for every $s,t\in \erre$.
We want to compare
the behaviour of the functions $|I_{\varepsilon,a}|^2$ and 
$\gamma^2$; to this aim, we define the function
$$\di (u,v,s,t) := |I_{\varepsilon,a}(u,v,s,t)|^2 - \gamma^2(u,v).$$
We have already shown (condition (e)) that 
\begin{equation}\label{zero}
\di(u,0,-u,u)=0.
\end{equation}
We start by proving that, if $V$ is sufficiently small, condition (d)
holds for every $(u,v)\in V$, for $t_1$ close to $-u$
and $t_2$ close to $u$. 
Using the definition of $\phi(u,v,z)$ on $A_1$ and $A_5$,
one can compute explicitly $\di(u,v,s,t)$
for $|s+u|\leq \alpha(-v)$ and for $|t-u|\leq \alpha(v)$.
By direct computations one obtains
\begin{equation}\label{deriv1}
\nabla_{v,s ,t}\, \di (u,0,-u,u) =0
\end{equation}
for $(u,0)\in V$.

We now want to compute the hessian matrix $\nabla^2_{v,s ,t}\, \di$ at the point
$(u_0, 0, -u_0, u_0)$.
By (\ref{sigma}) and (\ref{gi}), after some easy computations, we get
$$\diff_v^2\sigma (u,0)=
\frac{1}{2(u-a)}\diff_u \gamma(u,0)=
\frac{1}{2(u-a)}
\diff_v\diff_u \eta (u,0).$$
Using this equality and the explicit expression of 
$\di$ near $(u_0,0,-u_0,u_0)$, we obtain
\begin{eqnarray*}
\diff^2_v \di (u_0,0,-u_0,u_0) & = & -\frac{8\varepsilon}{(u_0-a)^2} 
(\diff_v \eta(u_0,0) - 4\varepsilon) +   \\
& & +\frac{2}{u_0-a} \diff_{v}\eta (u_0,0)\, \diff_v\diff_u \eta(u_0,0)
 - \diff^2_v (\gamma^2)(u_0,0).
\end{eqnarray*}
Since $\eta$ and
$\gamma$ do not depend on $a$ and 
$\varepsilon$,  for every $\varepsilon$
satisfying (\ref{gamma}) we can find 
$a$ so close to $u_{0}$ that 
\begin{equation}\label{det1}
\diff^2_v \di (u_0,0, -u_0, u_0) < 0.
\end{equation}
Moreover, we easily obtain that
\begin{eqnarray*}
&\displaystyle \diff_t^2 \di (u_0,0,-u_0,u_0) = 
\diff_s^2 \di (u_0,0,-u_0,u_0) =
8 - \frac{4}{\varepsilon}\diff_v \eta(u_0,0), & \\
&\displaystyle  \diff_v\diff_t \di (u_0,0,-u_0,u_0)  =
\diff_v\diff_s \di (u_0,0,-u_0,u_0) =
-\frac{4}{ u_0-a}(\diff_v \eta(u_0,0)- 4\varepsilon) , & \\
& \diff_t\diff_s \di (u_0,0,-u_0,u_0)  = 8. &
\end{eqnarray*}
By the above expressions,
it follows that 
\begin{eqnarray*}
\lefteqn{\det \left( \begin{array}{cc}
\diff_v^2 \di & \diff_v\diff_t \di  \\
& \\
\diff_v\diff_t \di  & \diff^2_t \di 
\end{array} \right)(u_0,0,-u_0,u_0) = } \\
 & = & \frac{16}{(u_0-a)^2}\diff_v \eta(u_0,0)
(\diff_v \eta(u_0,0) -4\varepsilon ) 
+\frac{c_1(\varepsilon)}{u_0-a} + c_2(\varepsilon),
\end{eqnarray*}
where $c_1(\varepsilon)$, $c_2(\varepsilon)$ 
are two constants depending only on
$\varepsilon$. Then, if $\varepsilon$ satisfies (\ref{gamma}), 
$a$ can be chosen so close to $u_0$ that
\begin{equation}\label{det2}
\det \left( \begin{array}{cc}
\diff_v^2 \di & \diff_v\diff_t \di  \\
& \\
\diff_v\diff_t \di  & \diff^2_t \di 
\end{array} \right)(u_0,0,-u_0,u_0)>0.
\end{equation}
At last, the determinant of the hessian matrix 
of $\di$ at $(u_0,0,-u_0,u_0)$ is given by
$$\det \nabla^2_{v,s ,t}\, \di (u_0,0,-u_0,u_0)= 
\frac{1}{u_0-a} (\diff_v \eta(u_0,0))^2
\diff_v \diff_u\eta(u_0,0)
(\diff_v \eta(u_0,0) - 4\varepsilon)
\frac{32}{\varepsilon^2 }  + c_3(\varepsilon),$$
where $c_3(\varepsilon)$
is a constant depending only on
$\varepsilon$. 
Since, by (\ref{psijacob}), 
$$\diff_v \diff_u\eta(u_0,0) = - \frac{\diff^2_x u(0,0)}{(\diff_x u(0,0))^3},$$
given $\varepsilon$ satisfying (\ref{gamma}), we can choose $a$ 
so close to $u_0$ that
\begin{equation}\label{det3}
\det \nabla^2_{v,s ,t}\, \di (u_0,0,-u_0,u_0)<0.
\end{equation}
By (\ref{det1}), (\ref{det2}), and (\ref{det3}),
we can conclude that, by a suitable choice of the parameters, the hessian matrix
of $\di$ (with respect to $v,s,t$) at $(u_0,0,-u_0,u_0)$ is negative definite. This fact, with 
(\ref{zero}) and (\ref{deriv1}), allows us to state 
the existence of a constant $\tau >0$
such that 
\begin{equation}\label{stima0}
\di(u,v,s,t)<0
\end{equation} 
for $|s+u_0|<\tau$, $|t-u_0|<\tau$,
$(u,v)\in V$, $v\neq 0$,
provided $V$ is sufficiently small.
So, condition (d) is satisfied for $|t_1+u_0|<\tau$ and $|t_2-u_0|<\tau$.
We can assume $\delta < \tau < \alpha (v)$ for every
$(u,v)\in V$.

From now on, since at this point the parameters $\varepsilon$, $a$
have been fixed, we simply write $I$ 
instead of $I_{\varepsilon, a}$. 
We now study the more general case $|t_1+u|<\alpha(-v)$ and $|t_2-u|<\alpha(v)$.

Let us set
$$m_1(u,v):= \max \left\{ | I(u,v,s,t ) | :
\, |s+u|\leq\alpha(-v),\, |t-u|\leq \alpha(v),\, |t-u_0|\geq\tau \right\}.$$
By the definition of $A_1, \ldots, A_5$, for
$\rho= \alpha(\delta) + \delta$ we have $(\phi^u,\phi^v)=0$ on
$(V { \times } [u_0 - \rho, u_0 +\rho])\setminus A_1$
and $(V { \times } [-u_0 - \rho,- u_0 +\rho])\setminus A_5$.
This implies that
$$m_1(u,v):= \max \left\{ | I(u,v,s,t ) | :
\, |s+u_0|\leq\rho,\, \tau \leq |t-u_0|\leq \rho \right\}$$
for $(u,v)\in V$.
The function $m_1$, as supremum of a family of continuous functions, is
lower semicontinuous. Moreover, $m_1$
is also upper semicontinuous; indeed,
suppose, by contradiction, that there exist two sequences
$(u_n)$, $(v_n)$ converging respectively to $u$, $v$,
such that $(m_1(u_n, v_n))$ converges to a limit
$l>m_1(u, v)$; then, there exist $(s_n)$, $(t_n)$
such that 
\begin{equation}\label{tnsn}
|s_n+u_n|\leq \alpha(-v_n), \qquad |t_n-u_n|\leq \alpha(v_n),
\qquad |t_n-u_0|\geq\tau,
\end{equation}
and $m_1(u_n, v_n)=| I(u_n,v_n,s_n,t_n ) |$.
Up to subsequences, we can assume that $(s_n)$, $(t_n)$
converge respectively to $s$, $t$ such that,
by (\ref{tnsn}), 
$$|s+u|\leq \alpha(-v), \qquad |t-u|\leq \alpha(v),
\qquad |t-u_0|\geq\tau;$$
hence, we have that
$$m_1(u, v)\geq | I(u,v,s,t ) |= \lim_{n\to\infty}
| I(u_n,v_n,s_n,t_n ) | =l>m_1(u, v),$$ 
which is impossible.
Therefore, $m_1$ is continuous.

Let $B$ be the open ball of radius $4\varepsilon$ 
centred at $(0,-4\varepsilon)$. Arguing as in (\ref{109}), we can prove that
\begin{equation}\label{901}
I(u,v,u,t) \in B
\end{equation}
whenever $0<|t-u|\leq\alpha(v)$.
In the same way we can prove that
\begin{equation}\label{902}
I(u,v,s,-u) \in B
\end{equation}
for $0<|s+u|\leq\alpha(-v)$.
We can write
\begin{equation}\label{spezza}
I(u,v,s,t) = I(u,v,s,-u) + I(u,v,-u,u) + I(u,v,u,t).
\end{equation}
So, for $|s+u|\leq\alpha(-v)$, $|t-u|\leq \alpha(v)$, and $|t-u_0|\geq\tau$,
by (\ref{902}), (\ref{phiu}), (\ref{phiv}), and (\ref{901}),
we obtain that
$$I(u,0,s,t)\in (0, \gamma(u,0)) + B + \overline{B},$$
hence, by (\ref{gamma}),
$I(u,0,s,t)$ belongs to the open ball of radius $\gamma(u,0)$
centred at $(0,0)$, and so, $m_1(u,0) < \gamma(u,0)$.
By continuity, if $V$ is small enough, 
\begin{equation}\label{emme1}
m_1(u,v) < \gamma(u,v)
\end{equation}
for every $(u,v)\in V$.

Analogously, we define
$$m_2(u,v):= \max \left\{ | I(u,v,s,t ) | :
\, |s+u|\leq\alpha(-v), |s+u_0|\geq\tau, |t-u|\leq \alpha(v),  \right\}.$$
Arguing as in the case of $m_1$, we can prove that, if
$V$ is small enough, 
\begin{equation}\label{emme2}
m_2(u,v) < \gamma(u,v)
\end{equation}
for every $(u,v)\in V$.

By (\ref{emme1}), (\ref{emme2}), and
(\ref{stima0}), we can conclude that
$I(u,v,t_1,t_2)$ belongs to the ball centred at $(0,0)$ with 
radius $\gamma(u,v)$, 
for $|t_1 +u|\leq \alpha(-v)$ and $|t_2 -u|\leq \alpha(v)$.
More precisely, let $E(u,v)$ be the intersection of this ball 
with the upper half plane 
bounded by the horizontal straight line passing through the point 
$(0, \frac{3}{4}\gamma(u,v))$:
by (\ref{spezza}), (\ref{stima2}), (\ref{901}), (\ref{902}), and    
(\ref{gamma}),
we deduce that  
\begin{equation}\label{911}
I(u,v,t_1,t_2)\in E(u,v)
\end{equation}
for $|t_1 +u|\leq \alpha(-v)$ and $|t_2 -u|\leq \alpha(v)$.

We can now conclude the proof of (d).
It is enough to consider the case $-u-\alpha(-v) \leq t_{1} \leq t_{2} 
\leq u + \alpha(v)$. We can write
\begin{eqnarray}
I(u,v,t_{1},t_{2}) & = &  I(u,v, t_{1} 
\land (-u+\alpha(-v)),
t_{2} \lor (u-\alpha(v))) + \nonumber \\
& & + I(u,v, t_{1} \lor (-u
+\alpha(-v)), t_{2} \land (u -\alpha(v)) )- \label{somma} \\
& & - I(u,v,-u +\alpha(-v), u-\alpha(v)). \nonumber
\end{eqnarray}

By (\ref{911}), it follows that  
\begin{equation}\label{inter}
I(u,v, t_{1} 
\land (-u+\alpha(-v)),
t_{2} \lor (u-\alpha(v)))\in E(u,v).
\end{equation}

Let $C_{1}(u,v)$ be the parallelogram having three consecutive vertices at 
the points 
$$(2hf(v), 0), \qquad (0,0), \qquad 
\sigma(u,v)\frac{(-v,u-a)}{\sqrt{(u-a)^{2}+v^{2}}},$$
let $C_{2}(u,v)$ be the segment with endpoints
$$(2hf(v),0), \qquad (2hf(v), 0) +
2\sigma(u,v)\frac{(-v,u-a)}{\sqrt{(u-a)^{2}+v^{2}}},$$
and let $C(u,v) :=C_{1}(u,v)\cup C_{2}(u,v)$.

From the definition of $\varphi$ in $A_{2}$, $A_{3}$, $A_{4}$, it 
follows that
\begin{equation}\label{922}
I(u,v, -u+\alpha(-v), u-\alpha(v)) = (2hf(v), 0) +
2\sigma(u,v)\frac{(-v,u-a)}{\sqrt{(u-a)^{2}+v^{2}}}
\end{equation}
and
\begin{equation}\label{923}
I(u,v, s_{1}, s_{2}) \in C(u,v)
\end{equation}
for $-u+\alpha(-v) \leq s_{1} \leq s_{2} \leq u- \alpha(v)$.
Let 
$$D(u,v) := C(u,v) - (2hf(v), 0) -
2\sigma(u,v)\frac{(-v,u-a)}{\sqrt{(u-a)^{2}+v^{2}}}.$$

From (\ref{somma}), (\ref{inter}),
(\ref{922}), and (\ref{923})
we obtain 
\begin{equation}
I(u,v , t_{1}, t_{2}) \in E(u,v) + D(u,v).
\end{equation}
As $|v| < \delta < 10\delta < u-a$ by (\ref{param}), 
the angle 
that the segment $C_{2}(u,v)$ forms with the vertical is less 
than $\arctan (1/10)$. 
Moreover, we may assume 
that the lenght $2\sigma(u,v)$ of the segment $C_{2}(u,v)$ is 
less than $\gamma(u,v)$; indeed, this is true for $v=0$ and, by 
continuity, it remains true if $\delta$ is small enough.
By (\ref{gamma}) and (\ref{stima4}), we have also
that $|2hf(v)| \leq \gamma(u,v) /16$.
Using these properties and simple geometric considerations, it is 
possible to prove that 
$E(u,v) + D(u,v)$
is contained in the ball with centre $(0,0)$ and radius $\gamma(u,v)$.
This concludes the proof of (d). 

If $\diff_x^2 u(0,0) <0$,
it is enough to change the definition of $\phi$ in the
sets $A_2$ and $A_4$, as follows:
$$\lambda \sigma(u,v) \frac{v}{\sqrt{(a-u)^2 + v^2}}\ddu
+ \lambda \sigma(u,v) \frac{a-u}{\sqrt{(a-u)^2 + v^2}}\ddv
+ \mu \ddz,$$
where $a> u_0 + 11\delta$ and
$$\sigma(u,v):= \frac{1}{2} \gamma(a-\sqrt{(a-u)^2 + v^2},0)-2\varepsilon.$$
\qed

\

\begin{theorem}\label{harm2}
Let $u: \Omega \to \erre$ be a harmonic function such that $\diff_y u(x,0) = 0$ for 
$(x,0) \in \Omega$, and let $w: \Omega \to \erre$ be the function defined by
$$w(x,y) := \begin{cases}
u(x,y) + 1 & \text{for $y>0$}, \\
u(x,y) & \text{for $y<0$}. 
\end{cases}$$
Assume that 
$\diff_x u(0,0)\neq 0$ and
$\diff_x^2 u(0,0)\neq 0$. 
Then there exists an open neighbourhood $U$ of $(0,0)$ such that $w$ is 
a Dirichlet minimizer in $U$
of the Mumford-Shah functional (\ref{g1}). 
\end{theorem}

{\bf Proof.}
We will write the 
calibration $\varphi$ as in (\ref{calibra2}) and we will adopt the representation 
(\ref{calibra}) for $\phi$.
We will use the same technique as in Theorem \ref{harm}. 
We give only the new definitions of the sets $A_1, \ldots, A_5$ and of the 
function $\phi$ when $\diff_x u(0,0)>0$ and $\diff^2_x u(0,0)>0$, 
and leave to the reader the verification of the fact that
this function is a calibration for suitable values of the involved
parameters. The case $\diff^2_x u(0,0)<0$ can be treated by the changes
introduced at the end of Theorem \ref{harm}. 

Let $u_0 := u(0,0)$.
Given $\varepsilon >0$, $h>0$, $\lambda >0$, and assuming 
$$V := \{ (u,v) : |u - u_0 | < \delta, |v| < \delta \},$$
we consider 
the following subsets of $V { \times } \erre$
\begin{eqnarray*}
A_1 & : = & \{ (u,v,z) \in V { \times } \erre : u + 1 - \alpha (v) < 
z < u + 1 + \alpha(v) \}, \\
A_2 & : = & \{ (u,v,z) \in V { \times } \erre : 5h + \beta(u,v) < z < 5h + \beta(u,v)
+ 1/ \lambda \}, \\
A_3 & : = & \{ (u,v,z) \in  V { \times } \erre : 2h < z < 4h \}, \\ 
A_4 & : = & \{ (u,v,z) \in  V { \times } \erre : h + \beta(u,v) < z < h + \beta(u,v)
+ 1/ \lambda \}, \\
A_5 & : = & \{ (u,v,z) \in  V { \times } \erre :  u  - \alpha (-v) <
z < u + \alpha(-v) \},
\end{eqnarray*}
where
$$\alpha(v) :=  \sqrt{4\varepsilon^2 - (\varepsilon -v)^2},$$
and $\beta$ is a suitable smooth function satisfying
$\beta (u,0) =0$, which will be defined later.
For $(u,v) \in V$ and $z\in\erre$ the vector
$\phi(u,v,z)$ is defined
as follows:
$$
\begin{cases}
\displaystyle \frac{2(\varepsilon- v)}{\sqrt{(\varepsilon-v)^2+(z-u-1)^2}} \ddu
-  \frac{2(z-u-1)}{\sqrt{(\varepsilon-v)^2
+(z-u-1)^2}} \ddv + \ddz
& \text{in $A_1$}, \\
\\
\displaystyle -\lambda \sigma (u,v) \frac{v}{\sqrt{(u-a)^2 + v^2}} \ddu +
\lambda \sigma (u,v) \frac{u-a}{\sqrt{(u-a)^2 + v^2}} \ddv +
\mu \ddz
& \text{in $A_2$}, \\
\\
\displaystyle  f(v)\ddu + \ddz  & \text{in
$A_3$}, \\
\\
\displaystyle -\lambda \sigma (u,v) \frac{v}{\sqrt{(u-a)^2 + v^2}} \ddu +
\lambda  \sigma (u,v) \frac{u-a}{\sqrt{(u-a)^2 + v^2}} \ddv +
\mu \ddz
& \text{in $A_4$}, \\
\\
\displaystyle  \frac{2(\varepsilon+ v)}{\sqrt{(\varepsilon+ v)^2 + (z- u)^2}} \ddu
+ \frac{2(z- u)}{\sqrt{
(\varepsilon+ v)^2 + (z- u)^2}} \ddv + \ddz  &
\text{in $A_5$}, \\
\\
\displaystyle \ddz  & \text{\hphantom{in $A_{5}$,}\llap{otherwise,}}
\end{cases}$$
where $a < u_{0}- 11 \delta$, $\mu >0$,
\begin{eqnarray*}
& \displaystyle f(v) := - \frac{1}{h} \left( \int_0^{\alpha(v) }
\frac{(\varepsilon - v)}{\sqrt{t^2 + (\varepsilon - v)^2}}\, dt +
\int_0^{\alpha(-v) } 
\frac{(\varepsilon + v)}{\sqrt{t^2 + (\varepsilon + v)^2}}\, dt
\right), &  \\
\\
\displaystyle & \sigma (u,v) := \frac{1}{2}
\gamma ( a + \sqrt{(u-a)^2 + v^2},0)
- 2 \varepsilon, & 
\end{eqnarray*}
and $\beta$ is the solution of the
Cauchy problem (\ref{cauchy}).
\qed

\end{document}